\documentstyle[12pt,righttag]{amsart} 
\input amssym.def		
\input amssym.tex 		

\newtheorem{prop}{Proposition}
\newtheorem{th}[prop]{Theorem}
\newtheorem{cor}[prop]{Corollary}
\newtheorem{lem}[prop]{Lemma}

\theoremstyle{definition}
\newtheorem{ack}{Acknowledgments} 
 
\newtheorem{qn}{Question} 
\theoremstyle{remark}
\newtheorem{rem}{Remark} 

\newtheorem{example}{Example}

\newcommand{\bbN}{{\Bbb{N}}}
\newcommand{\bbC}{{\Bbb{C}}}
\newcommand{\bbR}{{\Bbb{R}}}

\newcommand{\bbZ}{{\Bbb{Z}}}
\newcommand{\bbT}{{\Bbb{T}}}
\newcommand{\bbK}{{\Bbb{K}}}

\newcommand{\calW}{{\cal{W}}}
\def\real{{\Bbb R}}
\def\zed{{\Bbb Z}}

\newcommand{\be}{\beta}

\newcommand{\e}{\varepsilon}

\newcommand{\m}{\mu}

\newcommand{\r}{\varrho}
\newcommand{\s}{\sigma}
\newcommand{\f}{\varphi}

\newcommand{\Om}{\Omega}

\newcommand{\const}{\operatorname{const.}}

\newcommand{\Iso}{\operatorname{Iso}}

\newcommand{\disp}{\displaystyle}
\newcommand{\lb}{\label}

\newcommand{\lra}{\longrightarrow}

\newcommand{\wtw}{if and only if }

\newcommand{\Buo}{Without loss of generality }

\makeatletter
\def\@currentlabel{2.1}\label{e:dispaa}
\def\@currentlabel{2.21}\label{e:dispau}
\def\@currentlabel{2.22}\label{e:dispav}
 \def\@currentlabel{2.23}\label{e:dispaw}
\def\@currentlabel{2.24}\label{e:dispax}
\makeatother

\makeatletter
\def\alphenumi{%
  \def\theenumi{\alph{enumi}}%
  \def\p@enumi{\theenumi}%
  \def\labelenumi{(\@alph\c@enumi)}}
\makeatother

\begin{document}
\vspace{-1in}
\title{Isometric classification of norms  in rearrangement-invariant
function spaces}
\author{Beata Randrianantoanina}
\address{Department of Mathematics  \\ The University  of Texas at Austin
\\      Austin, TX 78712}
\email{brandri@@math.utexas.edu}
\maketitle
\begin{abstract}
Suppose that a real nonatomic function space on $[0,1]$ is equipped with two re\-arran\-ge\-ment-invariant norms $\|\cdot\|$ and $|||\cdot|||$. We study the question whether or not the fact that $(X,\|\cdot\|)$ is isometric to $(X,|||\cdot|||)$ implies that  $\|f\|= |||f|||$ for all $f$ in $X$. We show that in strictly monotone Orlicz and Lorentz spaces this is equivalent to asking whether or not the norms are defined by equal Orlicz functions, resp. Lorentz weights. 

We show that the above implication holds true in most rearrangement-invariant spaces, but we also identify a class of Orlicz spaces where it fails. We provide a complete description of Orlicz functions $\varphi \neq\psi$ with the property that $L_\varphi$ and $L_\psi$ are isometric.
\end{abstract}

\section{Introduction}


In this paper we study the following question:
\begin{qn}
Suppose that a   nonatomic function space on $[0,1]$ is equipped with two re\-arran\-ge\-ment-invariant norms $\|\cdot\|$ and $|||\cdot|||$ such that $(X,\|\cdot\|)$ and    $(X,|||\cdot|||)$ are isometric. Does it imply that the norms $\|\cdot\|$ and $|||\cdot|||$ are same? 
\end{qn}

Here  the word ``same'' could be understood in two ways:
\begin{itemize}
\item[$(a)$] we could say that $\|\cdot\|$ and $|||\cdot|||$ are same if  $\|f\|= |||f||| $ for all $f$ in $X$, i.e.  if the identity map $Id: (X,\|\cdot\|)\lra(X,|||\cdot|||)$ is an isometry, or
\item[$(b)$] if both $\|\cdot\|$ and $|||\cdot|||$ are Orlicz or Lorentz norms we could say that they are same if they are defined by  equal Orlicz functions or Lorentz weights, respectively.
\end{itemize}

It is well-known that if we do not require $|||\cdot|||$ to be rearrangement-invariant then the answer to Question~1 is no in either the sense $(a)$ or $(b)$, even when $(X,\|\cdot\|)= L_p[0,1]$  with the usual norm (see e.g. \cite{LW}).

Question~1$(a)$  has been previously studied with the additional assumption that $(X,\|\cdot\|) = (L_p[0,1], \|\cdot\|_p)$ (see \cite{AZ,AZ79,BH,HKW}).

Question ~1$(b)$   was asked by S. Dilworth and H. Hudzik.

Somewhat to author's surprise the answer to both Question~1$(a)$ and 1$(b)$ is negative -- there exist Orlicz functions $\f\neq\psi$ (then, clearly, $\|\cdot\|_\f, \|\cdot\|_\psi$ are different also in the sence $(a)$) so that $L_\f$ and $L_\psi$ are isometric.

In fact, Question~1$(a)$ and 1$(b)$ are equivalent for strictly monotone Orlicz-Lorentz spaces $X, \ Y$, i.e. $Id: X\lra Y$ is an isometry \wtw the Orlicz functions $\f_X,\ \f_Y$ and Lorentz weights $w_X, \ w_Y$ coincide (Theorem~\ref{LO}; this is not true in general, see \cite{HKW}). To our surprise the proof is much less obvious than one might expect.

Question~1$(a)$ for $L_p[0,1]$ has been studied by Abramovich and Zaidenberg \cite{AZ79,AZ}. They proved that if $Y$ is a (real or complex) rearrangement-invariant nonatomic function space on $[0,1]$ isometric to $L_p[0,1]$ for some $1\le p<\infty$ then the isometric isomorphism can be established via an identity map, i.e. $\|f\|_Y=\|f\|_p$ for all $f\in Y$ (cf. also \cite{HKW,BH} for the case when we additionally assume that $Y$ is an Orlicz space on $[0,1]$).

Zaidenberg \cite{Z77,Z80,Z95}  studied the general form of isometries between two complex rearrangement-invariant  spaces (r.i. spaces) and Jamison, Kami\' nska and P. K. Lin \cite{JKL} studied isometries between two complex Musielak-Orlicz spaces and between two real Nakano spaces. They proved that surjective isometries in such settings have to be weighted composition operators and Zaidenberg  characterized situations when the existence of isometry between complex r.i. spaces $X$ and $Y$ implies that the identity map is also an isometry.

Theorem~\ref{main} below generalizes these results to real spaces on $[0,1]$. We then use Zaidenberg's characterization of isometry groups of r.i. spaces to characterize when the identity map between real r.i. spaces is an isometry.
We also give a full description of the exceptional case of Orlicz spaces which can be isometric even when their Orlicz functions are different (Corollary~\ref{orlicz} provides the relation between the Orlicz functions that has to be satisfied in that case).

\begin{ack}  I wish to thank S. Dilworth and H. Hudzik for drawing my attention to this problem, S. Dilworth for many interesting discussions, Y. Abramovich for his interest in this work, and N. Carothers for his warm hospitality  during my visit at the Bowling Green State University, where this work was started.
\end{ack}
\section{Preliminaries}

Let us suppose that $\Omega$ is a Polish space and that
$\mu$ is a $\sigma$-finite Borel measure on $\Omega.$ We use
the term K\"othe space in the sense of \cite[p. 28]{LT2}.  Thus a
{\bf K\"othe function space} $X$ on $(\Omega,\mu)$ is a Banach
space of (equivalence classes of) locally integrable Borel
functions $f$ on $\Omega$ such that:  \newline (1) If
$|f|\le |g|$ a.e. and $g\in X$ then $f\in X$ with $\|f\|_X
\le \|g\|_X.$\newline (2) If $A$ is a Borel set of finite
measure then $\chi_A\in X.$

The {\bf K\"othe dual} of $X$ is denoted $X';$ thus $X'$ is the
K\"othe space of all $g$ such that $\int
|f||g|\,d\mu<\infty$ for every $f\in X$ equipped with the
norm ${\disp \|g\|_{X'}=  \sup_{\|f\|_X\le 1}\int |f||g|\,d\mu.}$
Then $X'$ can be regarded as a closed subspace of the dual
$X^*$ of $X$.   

A {\bf rearrangement-invariant function space} ({\bf r.i. space}) is a K\"othe function space on $(\Om,\mu)$   which satisfies the conditions:\newline
(1) $X'$ is a norming subspace
of $X^*.$
  \newline (2) If $\tau:\Om\lra \Om$ is any
measure-preserving invertible Borel automorphism then $f\in
X$ if and only if $f\circ\tau\in X$ and
$\|f\|_X=\|f\circ\tau\|_X.$ \newline (3)
$\|\chi_{B}\|_X=1$ if $\mu(B) = 1.$

The commonly studied r.i. spaces are classical Lebesgue spaces $L_p,$ Orlicz, Lorentz and Orlicz-Lorentz spaces. We recall the definitions below.

We say that $\varphi:\lbrack 0,\infty)\lra\lbrack 0,\infty)$\ is an
{\bf Orlicz function\/} if $\varphi$\ is non-decreasing and convex with
$\varphi(0)=0$. We define the {\bf Orlicz
space\/}
$L\sb \varphi$\ to be the space of those measurable functions $f$\ for
which $\|
f\|\sb \varphi$\ is finite, where $\|
f\|\sb \varphi$\  denotes the {\bf Luxemburg norm\/} defined by
$$  \|{f}\|\sb\varphi = \inf\left\{\, c :
   \int\sb \Omega \varphi\bigl(\frac{|{f(\omega)}|}{c}\bigr) \,d\mu(\omega) \le 1
\,\right\}.$$

 If $f$\ is a
measurable function, we define the
{\bf non-increasing rearrangement\/} of $f$\ to be
$$ f\sp *(x) = \sup\bigl\{\, t : \mu(| f| \ge t) \ge x \,\bigr\}.$$
If $1\le q <
\infty,$  
and if $w:(0,\infty)\to(0,\infty)$\ is a non-increasing
function, we
define the {\bf Lorentz space\/} $L\sb {w,q}$\ to
be the space
of those measurable functions $f$\ for which $\|f\|\sb {w,q}$\ is
finite,  where $\|f\|\sb {w,q}$ denotes the 
{\bf Lorentz norm\/} defined by
$$ \| f\|\sb {w,q} = \left(\int\sb 0\sp \infty w(x) f\sp *(x)\sp q
\,dx\right)\sp {1/q}.$$

  If $\varphi$\
is an Orlicz function, and if
$w:\lbrack 0,\infty)\to\lbrack 0,\infty)$\ is a non-increasing
function, we define the {\bf Orlicz--Lorentz space\/} $\Lambda\sb
{w,\varphi}$\ to be
the   space of measurable functions $f$\ for which $\| f\|\sb
{w,\varphi}$\ is
finite, where $\| f\|\sb{w,\varphi}$\ denotes 
 the {\bf Orlicz--Lorentz norm\/} defined by $$ \| f\|\sb {w,\varphi} = \inf\left\{\, c :
   \int\sb 0\sp \infty w(x) \varphi\bigl(\frac{f\sp *(x)}{c}\bigr) \,dx \le 1
\,\right\}.$$

 The {  Orlicz--Lorentz
spaces} are a common generalization of the Orlicz spaces and the Lorentz
spaces.

An operator $T:X\lra X$ will be called {\bf elementary} or a {\bf weighted composition operator} \lb{el} if
there is a Borel function $h$ and a Borel map
$\sigma:\Omega\lra\Omega$ such that $Tf(s)=h(s)f(\sigma(s))$
a.e. for every $f\in X.$ Observe that a necessary condition
on $a$ and $\sigma$ is that if $B$ is a Borel set with
$\mu(B)=0$ then $\mu(\sigma^{-1}B\cap\{|h|>0\})=0.$

We now need to introduce a technical definition.  We will
say that an r.i. space $X$ has {\bf property ${\bold(GP)}$} if there exists $n\in\bbN$ such that for every
$t>0,$ $$ \|\chi_{[0,2^{-n}]}\|_X < \|\chi_{[0,2^{-n}]}+t\chi_{[2^{-n},1]}\|_X.$$ We say that
$X$ has property ${\bold(GP')}$ if $X'$ has property (GP).

Notice that if $X$ is strictly monotone then $X$ has   property $ (GP)$. The reason for introducing this property, rather than simply dealing with strictly monotone spaces, is that every rearrangement-invariant space $X$ has to satisfy at least one of the properties $(GP)$ or $(GP')$. Indeed, if $X$ fails both $(GP)$ and $(GP')$, then for all $n\in\bbN$, say $n=1$, there exists $\eta>0$ small enough so that 
$ \|\chi_{[0,2^{-1}]}\|_X =\|\chi_{[0,2^{-1}]}+\eta\chi_{[2^{-1},1]}\|_X$
and $ \|\chi_{[0,2^{-1}]}\|_{X^*} =\|\chi_{[0,2^{-1}]}+\eta\chi_{[2^{-1},1]}\|_{X^*}$. But then
$$\frac{1}{2}(1+\eta^2) = \int(\chi_{[0,2^{-1}]}(s)+\eta\chi_{[2^{-1},1]}(s))^2\ ds \le \|\chi_{[0,2^{-1}]}\|_X \|\chi_{[0,2^{-1}]}\|_{X^*}= \frac{1}{2}$$
which contradicts the fact that $\eta>0$.

Notice that for any $p$, $1<p<\infty,$ $L_p$ satisfies both $(GP)$ and $(GP')$.

An Orlicz space $L_\varphi$ satisfies  $(GP)$  \wtw $\varphi(t)>0$ for all $t>0$ and there exists $t>1$ with $\varphi(t)<\infty$.

A Lorentz space $L_{w,p},\ \ 1\le p <\infty,$ satisfies  $(GP)$  whenever there exists $t>0$ with $w(t)>0,$ i.e. whenever $L_{w,p}\neq L_\infty.$

An Orlicz--Lorentz
space $\Lambda_{w,\varphi}$ satisfies  $(GP)$  whenever $\varphi(t)>0$ for all $t>0$, there exists $t>1$ with $\varphi(t)<\infty$ and there exists $s>0$ with $w(s)>0.$

\section{Isometries between two different spaces $X$ and $Y$}\lb{s1}

We start with the statement of the theorem about the form of surjective
isometries between two different rearrangement-invariant spaces $X$ and $Y$.
The proof is a technical refinement of    \cite[Theorem 6.4]{KR}
and we present it in section~\ref{last}.

\begin{th} \lb{main}
Let $X, \ Y $ be nonatomic rearrangement-invariant spaces on $[0,1]$
which are not isometrically equal to $L_2[0,1],$ and
such that both $X$ and $Y$ satisfy property $(GP)$  or (both) $(GP').$
Suppose that
$T:X\lra Y$ is a surjective isometry.  Then there exists a
Borel function $h$ on $[0,1]$ with $|h|>0$ and an invertible
Borel map $\sigma:[0,1]\lra [0,1]$ such that
$\lambda(\sigma^{-1}(B))>0$ if and only if $\lambda(B)>0$
for $B\in\cal B$ and so that $Tf(s) = h(s)f(\sigma(s))\ \  a.e.$
for every $f\in X.$ \end{th}

\begin{rem}
Unfortunately we were not able to eliminate a technical requirement for
spaces $X,Y$ to satisfy property $(GP)$ or $(GP')$.
On one hand, as dicussed above, this is not a very restrictive assumption (and this is the  reason presenting it here in this
slightly unfinished form).
On the other hand there is a general question about the group 
$\calW$ of all invertible
weighted composition operations which is interesting in itself and which
would immediately imply Theorem~\ref{main}.
Namely if $T$ is an invertible operator such that $T\circ \calW\circ T^{-1}
\subset \calW$ does it imply that $T\in \calW$?
Our proof shows that the answer to this question is yes if we assume a special
form of $T$.
\end{rem}

Now we will apply the analysis of isometry groups of rearrangement-invariant
function spaces due to Zaidenberg \cite{Z77}.

Following his notation we denote by $NS= NS[0,1]$ the group of all invertible
mappings of $\s : [0,1]\to [0,1]$ which, together with their inverses, preserve
measurability.
All surjective isometries of $L_p$ ($p\neq2$) have form described in Theorem~\ref{main} (\cite{Lam}). So there is a 1-1 correspondence between the group of all (positive) surjective isometries of $L_p$ and the group $NS$. Let $NS_p$ denote the group $NS$ equipped with the  topology induced by weak convergence of operators in $L_p$.

If $\s : [0,1]\to [0,1]$ preserves measurability then we denote by
$\mu_\s$ the measure defined by $\mu_\s (A) = \mu(\s^{-1}(A))$.
Since $\mu_\s (A)=0$ iff $\mu(A)=0$ there exists a Radon-Nikodym
derivative $d\mu_\s\over d\mu$ of $\mu_\s$, and we will denote
$d\mu_\s\over d\mu$ by $\s'$.
If $G$ is a subgroup of $NS$ we denote by $N(G)$ the normalizer of 
$G$ in $NS$.
Finally if $X$ is a function space we denote by $\Iso (X)$ the group 
of invertible isometries of $X$. 

The following is due to Zaidenberg
(\cite[Theorems~2 and 3]{Z77}) (for proof see \cite{Lin95}).

\begin{th}\lb{Z}  \mbox{   }

\begin{enumerate}
\item[a)] If $X$ is a rearrangement-invariant function space on $[0,1]$
and $X\ne L_2[0,1]$ then the group of invertible isometries of $X$ coincides
with one of the following closed subgroups of $NS = NS_p$:
\begin{enumerate}
\item[(i)] $NS$
\item[(ii)] $NS (\real_+;a) =\{ \s \in NS \mid \exists b>0;
\s' (t) \in \{ba^k\}_{k=-\infty}^\infty$ for almost all $t\in [0,1]\}$,
where $a>1$.
\item[(iii)] $NS (a;d) = \{\s\in NS\mid \exists s\in\zed: \s' (t)
\in \{a^{s+kd}\}_{k=-\infty}^\infty$ for almost all $t\in [0,1]\}$,
where $a>1$, $d\in\zed$, $d\ge1$.
\item[(iv)] $U= \{\s \in NS\mid\s$ is measure preserving$\}$. 
\end{enumerate}
\item[b)] Two subgroups of $NS$ from the list (i)--(iv) are conjugate 
if and only if they coincide.
\item[c)] If $G$ is a subgroup of type (i), (ii) or (iv) then $N(G)=G$,
and $N(NS(a;d)) = NS(\real_+;a^d)$.
\end{enumerate}
\end{th}

Theorem~\ref{Z} allows us to precisely identify which weighted composition operators can act as surjective isometries between r.i. spaces.

\begin{prop}\lb{isolp}
Suppose that $X,Y$ are rearrangement invariant function spaces and 
$T:X \buildrel \mbox{\em onto}\over \longrightarrow Y$ is an isometry 
such that $Tf= hf\circ\sigma$. 
Then either $\sigma$ is measure preserving and
$|h| = 1$ a.e. or there exists $ p$, $1\le p\le\infty$,   
such that $X,Y$ are equal to $L_p$ with
some equivalent norms and $T$ is also an isometry as an operator from
$L_p$ to $L_p$, i.e. $|h(t)|^p = \s'(t) \ \ a.e.$
\end{prop}

\begin{pf}
We consider two cases. 
First we assume that $\Iso (X)=U$  (or, symmetrically, $\Iso (Y)=U$). 
If $\tau\in NS$ we denote by $V_\tau$ the composition map defined 
on $Y$ by 
$$V_\tau f= f\circ \tau\ .$$ 
If $\tau$ is measure preserving then $T^{-1}V_\tau T$ is an isometry of $X$. 
Hence we have 
$$T^{-1} V_\tau Tf(t) = {h\circ\tau (t) \over h\circ\sigma^{-1}(t)} 
\ f\circ \sigma \circ\tau\circ \sigma^{-1} (t)$$ 
and $\sigma \circ\tau\circ\sigma^{-1} \in U$ and $|({h\circ\tau(t))/(h\circ\sigma^{-1}(t))}| =1$ a.e. 
Since $\tau \in U$ is arbitrary we conclude that $|h| \equiv \const$ a.e. 
so $|h|\equiv 1$ a.e. and by Theorem~\ref{Z}(c) $\sigma \in U$, i.e., 
$\sigma$ is measure preserving. 

If $\Iso (X) \ne U$ then by    \cite[Theorem~7.2]{KR} (cf.\ also 
\cite{Z77}) $X$ equals to $L_p$ for some $p$, $1\le p\le\infty$, with an 
equivalent renorming and Boyd indices $p_X,q_X$ of $X$ are both equal to $p$. 
Similarly, $\Iso (Y)\ne U$ implies that $Y\cong L_p$ and $p_Y = q_Y =p$ 
with clearly the same $p$, since $T$ is an isometry. 

Now, by   \cite[ Proposition~7.1]{KR} (cf.\ also    \cite[Theorem~5.1]{KR}) 
$T$ is an isometry of $L_p$. 
We should remark that   \cite[Proposition~7.1]{KR} is formulated for 
isometries acting on only one space $X\to X$ but in fact it uses only the 
equality of Boyd indices of the range  and domain spaces. 
Since there is no change in the proof we will not repeat it here. 
\end{pf}

As an immediate consequence of Theorem~\ref{main} and Proposition~\ref{isolp} 
 we obtain the following corollary.

\begin{cor}\lb{az}
Suppose that $X, \ Y$ are nonatomic rearrangement-invariant spaces on $[0,1]$ not isometrically equal to $L_2[0,1],$ and
both satisfying property $(GP)$ or $(GP')$ 
and that there exists an isometry $T:X \buildrel \mbox{\em onto}\over \longrightarrow Y.$
Then
$$\Iso(X) = \Iso(Y).$$
Further, if $\text{Iso}(X) \neq NS(a,d)$ then
$Id: X\buildrel \mbox{\em onto}\over \lra Y $ is also
an isometry.
\end{cor}

In particular identity map is an isometry if one of the
spaces $X$ or $Y$ is equal to $L_p$ for some $p$, $1<p<\infty, p \ne 2,$ since 
then $L_p$ satisfies both $(GP)$ and $(GP')$. 
Thus we  generalize   the result of Abramovich and Zaidenberg
\cite{AZ} who proved   Corollary~\ref{az}  in the case when one of the spaces
$X$ or $Y$ equals to $L_p$, for some $1\le p\le\infty$.

Also notice that, by    \cite[Theorem~7.2]{KR}, if $X,Y$ are not isomorphic to
any $L_p$ then $\Iso (X) = \Iso (Y) =U$ and hence Corollary~\ref{az} 
holds.

Furthermore, since $N(NS(a;d)\supsetneq NS(a;d)$
(Theorem~\ref{Z}(c)) there are isometric spaces $X,Y$ such that
$\Iso (X) = \Iso(Y) = NS(a;d)$ and the identity map is not an isometry
--- notice that identity map is well defined since $X,Y$ are both equal
to $L_p$ with an equivalent norm.
Thus, as a consequence of Theorem~\ref{Z}(c), we get a positive answer
to Question~1(a). In fact, we even can get two Orlicz spaces with this property.

\begin{example} (see \cite{Lam,Z77} and Corollary~\ref{orlicz} below)

Consider the Orlicz function $\f(t) = t^5 \exp (\sin \ln t)$.
Take any $\sigma \in NS (\real_+,e^{2\pi\cdot 5})\setminus NS(e^{2\pi\cdot5};1)$
and let $h$ be such that
$Tf= h\ f\circ\sigma$ defines the isometry on $L_5$, i.e., $|h(t)|^5=
\sigma'(t)$ for almost all $t$.
Now define a new norm $\|\cdot\|_X$ on $L_5$ by $\|f\|_X = \|Tf\|_\f$.
Then, clearly $X$ and $L_\f$ are isometric.
Moreover 
\begin{align*}
\|f\|_X &= \| h\cdot f\circ \sigma\|_\f = \inf\left\{c\ :\ \int \f\left(\frac{|h(t) f(\s(t))|}{c}\right) \ dt \le 1\right\} \\
      &= \inf\left\{c\ :\ \int \frac{|h(t)|^5 |f(\s(t))|^5}{c^5}\cdot\exp\left[\sin\left(\ln |h(t)| + \ln\frac{|f(\s(t))|}{c} \right)\right] \ dt \le 1 \right\} \\
\intertext{and letting $u=\s(t)$, since $|h(t)|^5=
\sigma'(t)$}
&=\inf\left\{c\ :\ \int \frac{  |f(u)|^5}{c^5}\cdot \exp\left[\sin\left(\frac15\ln(\s'(t))+ \ln\frac{|f(u)|}{c} \right)\right] \ dt \le 1 \right\} \\
\intertext{since for all $t$ there exists $k(t)\in\bbZ$ with $\s'(t) = b e^{2\pi\cdot 5 k(t)}$ }
&= \inf\left\{c\ :\ \int \frac{  |f(u)|^5}{c^5}\cdot \exp\left[\sin\left(\ln(b^{\frac15}\cdot \frac{|f(u)|}{c} )\right) \right] \ dt \le 1 \right\}
\end{align*}

Thus $\|\cdot\|_X$ is   an Orlicz norm defined by the
Orlicz function
$\f_\sigma (t) = t^5 \exp (\sin\ln b^{1/5}t)$, where
$b\in \real_+ \setminus \{ e^{10\pi k}; k\in\zed\}$ is such that
$\sigma' (t) \in \{b e^{2\pi\cdot 5 k}:k\in\zed\}$.
\end{example}

\section{Isometric uniqueness of norm}\lb{s2}

In this section we study Question~1(b) for Orlicz-Lorentz spaces.
To our surprise the proof, while definitely not difficult, is not completely
obvious.

We will work not with the usual definition of Orlicz-Lorentz spaces (which was presented in the preliminaries)
but with the equivalent one introduced by Montgomery-Smith \cite{SMS}.
We feel that this definition is very close to the spirit of
Krasnosielskii-Rutickii's definition of Orlicz spaces (\cite{Kr-Rut}) and it is
convenient for us for technical reasons.
We briefly remind  the notations from \cite{SMS} below.

First we define $\varphi-$functions. These replace the notion of Orlicz
functions in
our discussions.

 A {\bf $\varphi-$function} is a continuous, strictly increasing function $F:[0,\infty)
\lra
[0,\infty)$ such that
 $F(0) = 0$; $F(1)=1$ 
and $\lim\sb {n\to\infty} F(t) = \infty$;

The definition of a $\varphi-$function is slightly more restrictive than
that of an
Orlicz function in that we insist that $F$\ be strictly increasing.

If $F$\ is a $\varphi-$function, we will define the function $\widetilde F(t)$\
to
be $1/F(1/t)$\ if $t>0$, and $0$\ if $t=0$.

 If $(\Omega,\cal F,\mu)$\ is a measure space, and
$F$\ and $G$\ are $\varphi-$functions, then we define the {\bf Orlicz--Lorentz functional} of a measurable function $f$ by
$$ \|{f}\|\sb {F,G} = \|{f\sp *\circ\tilde F\circ\widetilde G\sp
{-1}}\|\sb G .$$

The {\bf  Orlicz--Lorentz space} 
$L\sb {F,G}(\Omega,\cal F,\mu)$\ (or $L\sb {F,G}$\ for
short) is the vector space of measurable functions $f$\ for which
$\|
f\|\sb {F,G} <\infty$, modulo functions that are zero almost everywhere.

Notice that if $\Omega = [0,1]$ then functions $F$ and $G$ do not need to 
be defined on $[0,\infty)$. 
We require only that  $(\widetilde F\circ 
\widetilde G^{-1})^{-1} $ is defined on $[0,1]$, and thus $F$ and $G$ are really defined 
on $[0,1]$. 

The Orlicz--Lorentz spaces defined here are equivalent to the
definition given in the preliminaries. Namely, 
 if $w$\ is a weight function, and $G$\ is a $\varphi-$function, then
$\Lambda\sb {w,G} = L\sb {\tilde W\sp {-1} \circ G,G}$, where
$$ W(t) = \int\sb 0\sp t w(s) \, ds .$$

We will also need the following two simple lemmas from \cite{SMS}.
For the sake of completeness we present their short proofs with adjustments for our 
situation.

\begin{lem} \lb{5.2.1} Suppose that $F$, $G\sb 1$\ and $G\sb 2$\ are
$\varphi-$functions.
Then $\| f\|\sb {F,G\sb 1}
= \|
f\|\sb {F,G\sb 2}$\ for all measurable $f$\  if and only if
$\| f\|\sb {1,G\sb 1} = \| f\|\sb {1,G\sb 2}$\ for all
measurable
$f.$
\end{lem}

\begin{pf}  This follows because $\|f\|\sb {F,G} = \|{f\sp *\circ
\tilde
F}\|\sb {1,G}$.
\end{pf}

\begin{lem} \lb{5.2.2} Suppose that $G\sb 1$, $G\sb 2$\ and $H$\ are
$\varphi-$functions.
Then if $\| f\|\sb {1,G\sb 1} = \| f\|\sb {1,G\sb 2}$\ for all
measurable $f$, then $\| f\|\sb {1,G\sb 1\circ H} =\| f\|\sb {1,G\sb 2\circ
H}$\ for all measurable
$f.$\ \end{lem}

\begin{pf}  Suppose that $\| f\|\sb
{1,G\sb 1} =\| f\|\sb {1,G\sb 2}$\ for all measurable $f$. Fix $\e >0$ 
and let $g$ be any simple function such 
that $1/(1+\e)\le \| g\|\sb {1,G\sb 2\circ H} < 1$ and 
$g^* = \sum_{i=1}^n a_i \chi_{A_i}$. 

Then there exists $\lambda<1$ such that 
\begin{equation*}
\begin{split}
1&\ge \int G_2 \circ H\circ {1\over\lambda} g^* \circ \widetilde H^{-1} 
\circ \widetilde G_2^{-1} (x)\,dx 
 = \sum_{i=1}^n \int G_2\circ H\left( {1\over\lambda} a_i 
\chi_{\widetilde G_2\circ \widetilde H (A_i)} (x)\right)\,dx\cr 
&= \sum_{i=1}^n G_2 \circ H \left({1\over\lambda} a_i\right) 
\mu (\widetilde G_2 \circ\widetilde H(A_i))  
> \sum_{i=1}^n G_2 \circ H (a_i)\mu (\widetilde G_2\circ\widetilde H(A_i))\cr
&= \int G_2 \circ H\circ g^* \circ \widetilde H^{-1}\circ 
\widetilde G_2^{-1} (x)\,dx, 
\end{split}
\end{equation*} 
since $H$ is strictly increasing.

Therefore, $\|{H\circ g\sp * \circ \widetilde H\|\sp {-1}}\sb {1,G\sb 2}
< 1$, and hence
$\|{H\circ g\sp *\circ \widetilde H\sp {-1}}\|\sb {1,G\sb 1}< 1$. 
Therefore,
$$\int G\sb 1\bigl(H \circ g\sp *
   \circ \widetilde H\sp {-1} \circ \widetilde G\sb 1\sp {-1} (x)\bigr) \, dx
\le 1 .$$
that is, we have
$$ \int G\sb 1\circ H\bigl( g\sp *
   \circ \tilde H\sp {-1} \circ \widetilde G\sb 1\sp {-1} (x)\bigr) \, dx
\le 1 ,$$
that is, $\| g\|\sb {1,G\sb 1\circ H} \le 1.$
Hence $\| g\|\sb {1,G\sb 1\circ H} \le (1 + \e )\| g\|\sb {1,G\sb 2\circ H}.$ Since $\e$ is arbitrary and by symmetry we get that
$\| g\|\sb {1,G\sb 1\circ H} = \| g\|\sb {1,G\sb 2\circ H}$ for all simple functions.
\end{pf}

We are now ready to prove the main result of this section.

\begin{th}\lb{LO}
Suppose that $L_{F_1,G_1}, L_{F_2,G_2}$ are Orlicz-Lorentz spaces and
that the identity map is an isometry, i.e. for all $f\in L_{F_1,G_1}$,
$\|f\|_{F_1,G_1} = \|f\|_{F_2,G_2}$.
Then $F_1(t) = F_2(t)$  and $G_1(t) = G_2(t)$ for all $t$.
\end{th}

\begin{rem}
Notice that, by Corollary~\ref{az}, whenever $L_{F_1,G_1}, L_{F_2,G_2}$ are isometric and $\Iso(L_{F_1,G_1}) \neq N(a,d),$ for any $a>1, \ d\in\bbN,$ then
the identity map is an isometry. In particular this holds whenever $L_{F_1,G_1}, L_{F_2,G_2}$ are not isomorphic to any $L_p,$ and also when $L_{F_1,G_1}=L_p$.
\end{rem}
  
\begin{pf}
Let $A\subset \Omega$. Then for $i=1,2$
$$\|\chi_A\|_{F_i,G_i}
= \widetilde F_i (\mu (A) = {1\over F_i ((\mu(A))^{-1})}\ .$$
Hence if $\Omega = [0,\infty)$ then $F_1 (t) = F_2(t)$ for all $t\ge0$ 
and if $\Omega = [0,1]$ then $F_1(t)=F_2(t)$ for $t\ge 1$, 
but as we mentioned in the definition of $L_{F,G}([0,1])$, $[1,\infty)$ 
is the required domain of $F$ in this case. 
Now by Lemma~\ref{5.2.1} we get that for all $f$
$\|f\|_{1,G_1} = \|f\|_{1,G_2}$ and by Lemma~\ref{5.2.2}
$\|f\|_{1,G_1\circ G_2^{-1}} = \|f\|_{1,1}$.
So it is enough to prove that if $\|f\|_{1,G} = \|f\|_1$ for all $f\in L_1$
then $G(t) =t$ for all $t$.
For any $s$, $0<s<1$ and $a$, $0<a\le 1< {1\over s}$ consider
$f= {1-as\over 1-s} \chi_{[0,1-s]} + a\chi_{[1-s,1]}$.
Then $\|f\|_1 = 1$ and $f^*=f$.
Hence $\|f\|_{1,G} =1$.
Since
$$f\circ \widetilde G^{-1} (t) = {1-a+as\over s} \chi_{[0,{1\over G(1/s)}]}
+ a\chi_{[{1\over G(1/s)},1]}$$
and since $G$ is continuous, $\|f\|_{1,G}=1$ implies that for all
$s$, $0<s<1$ and $a$, $0<a\le 1< {1\over s}$:
\begin{equation}\lb{LO1}
{1\over G(1/s)} G\left( {1-a+as\over s}\right)
+ \left( 1-{1\over G(1/s)}\right) G(a) =1
\end{equation}
Now, for a fixed $s$, take a left derivative of equation~\eqref{LO1} with
respect to $a$.
We get
$${1\over G(1/s)} G' \left( {1-a+as\over s}\right)
\left( {-1+s\over s}\right) + \left( 1-{1\over G(1/s)}\right) G' (a)=0$$
where $G'(t)$ denotes the left derivative of $G$.
In particular when $a=1$ we have
$${1\over G(1/s)} G'(1) \left( {-1+s\over s}\right)
+ \left( 1-{1\over G(1/s)}\right) G'(1)=0\ .$$
Since $G$ is convex and strictly increasing $G'(1)\ne 0$, so
$${1\over G(1/s)} \left( {-1+s\over s}\right)
+ \left( 1-{1\over G(1/s)}\right) =0$$
and, after simplification
$G(1/s)=1/s$ for all $s$, $0<s<1$, i.e., $G(t)=t$ for all $t>1$ 
which ends the proof in the case when $\Omega = [0,1]$ since in that 
case $G$ is considered as a function from $[1,\infty)$. 

If $\Omega = [0,\infty)$ we consider consider any $a<1$. 
Then ${1-a+as\over s} >1$ and \eqref{LO1} can be rewritten as:
$$s\cdot {1-a+as\over s} + (1-s) G(a)=1\ .$$
and therefore
\begin{gather*}
1-a+as +(1-s)G(a)=1\\	
(1-s) G(a) = a(1-s)\\	
G(a) =a\mbox{ for all } a<1\ .	
\end{gather*}
By continuity $G(t)=t$ for all $t>0$.
\end{pf}

As a corollary of Theorem~\ref{LO} we obtain the characterizations of isometric Lorentz spaces and Orlicz spaces.

\begin{cor}\lb{lorentz}
Two Lorentz spaces $L_{w_1,p_1}[0,1], \ L_{w_2,p_2}[0,1]$ are isometric
\wtw $w_1(t)=w_2(t)$ a.e. and $p_1=  p_2$, or $L_{w_1,p_1}, \ L_{w_2,p_2}$ are isometric to $L_\infty$, i.e. for $i=1, 2$, $p_i = \infty$ or $w_i(t) = 0$ for all $t>0$.
\end{cor}

\begin{pf}
If there exist $t_1, t_2$ such that $w_i(t_i) >0$ for $i=1, 2$, and $p_1, p_2 \neq\infty$, then  $L_{w_1,p_1}, \ L_{w_2,p_2}$ satisfy property $(GP)$. Also, by \cite{CHL}, $\Iso(L_{w_1,p_1}) = \Iso(L_{w_2,p_2}) =U.$ Hence, by Corollary~\ref{az}, $Id:L_{w_1,p_1}\lra L_{w_2,p_2}$ is an isometry.

If $w_1(t),w_2(t)>0$ for all $t>0$, we can apply Theorem~\ref{LO} to get the conclusion. 
Otherwise let $$F_{s,i}(a) = \|\chi_{[0,s]} + a \chi_{[s,1]}\|_{w_i,p_i}$$
for $i=1, 2$. Since the identity map is an isometry $F_{s,1}(a) = F_{s,2}(a)$ for all $s, a \le 1.$ In particular
\begin{equation} \lb{lor1}
\left(\int_{[0,s]} w_1(t)\ dt\right)^{\frac{1}{p_1}}= F_{s,1}(0) = F_{s,2}(0) =\left(\int_{[0,s]} w_2(t) \ dt\right)^{\frac{1}{p_2}}.
\end{equation}
Moreover
\begin{equation} \lb{lor2}
 \int_{[s,1]} w_1(t) \ dt = F_{s,1}'(1) = F_{s,2}'(1) = \int_{[s,1]}w_2(t)\ dt.
\end{equation}

Since \eqref{lor1} and \eqref{lor2} hold for all $s\le 1$ we conclude that $p_1=  p_2$ and $w_1(t)=w_2(t)$ a.e.

If  $w_1(t) = 0$ for all $t>0$ or if $p_1 = \infty$, then $L_{w_1,p_1} = L_\infty$ and, by \cite{AZ}, $Id:L_{w_1,p_1}= L_\infty\lra L_{w_2,p_2}$ is an isometry, and thus also $p_2 = \infty$ or $w_2(t) = 0$ for all $t>0$.
\end{pf}

\begin{cor}\lb{orlicz}
Let ${\varphi}, {\psi}$ be two  Orlicz functions such that $0<\varphi(t)<\infty$ and $0<\psi(t)<\infty$ for all $t>0.$ Then Orlicz spaces $L_{\varphi}$ and $L_{\psi}$ are isometric \wtw one of the two possibilities holds:

\begin{enumerate}

\item[(1)] $\varphi (t) = \psi (t)$ for all  $t$.
\item[(2)] there exist $b>0$ and $p\ge1,$ such that $\f(t) = (1/b)\psi (b^{(1/p)}t) $
for all $t$.
In this case $L_\varphi,L_\psi$ are isomorphic to  $L_p$. 
\end{enumerate}
\end{cor}

\begin{rem}
If $\varphi (t)= \psi (t) =\infty$ for all $t>1$, then
  $L_\varphi$ and $L_\psi$ are  isometric with $L_\infty$ via an identity
map (\cite[Corollary 1]{HKW}).
\end{rem}

For the proof of the above corollary we need a description of isometry groups of Orlicz spaces, which is due to Zaidenberg \cite{Z76,Z77}. However, since the English translation of \cite{Z77} has an unfortunate misprint (the author has not seen the original Russian version) and since the proof is not presented in full there, we will enclose both the exact statement and its proof (we use the same notation as in Theorem~\ref{Z}).

\begin{lem} \lb{gror}
If $L_\f$ is a (real or complex) Orlicz space ($L_\f\neq L_2$) then the group $G_\f$ of all surjective isometries of $L_\f$ equals:
\begin{enumerate}

\item[(a)] $NS$, when $\f(t) = t^p$ for some $p\ge 1, \ p\neq2.$
\item[(b)] $NS(\bar a^p,1)$, when the multiplier group of $\f$ is generated by $\bar a, \ \bar a>1,$ and $L_\f $ is isomorphic to $L_p,$
\item[(c)]  $U$, otherwise.
\end{enumerate}
\end{lem}

\begin{pf*}{Proof of Lemma~\ref{gror}}
By \cite[Theorem~4]{Z76} (cf. also \cite{KR}) every surjective isometry $S$ of $L_\f$ has form
$$Sf = h \ f\circ \s,$$
where $\s$  an invertible Borel map $\s:[0,1]\lra[0,1]$ and $h:[0,1]\lra\bbK$ ($\bbK = \bbR$ or $\bbC$) are such that 
\begin{equation} \lb{fi}
\f(|h(t)|\cdot x) = \s'(t)\f(x) \text{\ \ \ \ for $a.e.\ t$ and all $x\ge0$}
\end{equation}

By Theorem~\ref{Z}$(a)$ the isometry group $G_\f$ coincides with one of the following groups:
$NS$, $NS(\bbR_+,a)$, $NS(a,1)$, $U.$

As mentioned above, it is well known that $G_\f = NS$ \wtw $L_\f = L_p.$

If $G_\f \neq U$ then $L_\f$ is isomorphic to $L_p$ (\cite[Theorem~4]{Z77}, \cite[Theorem~7.2]{KR}). Moreover, in this case we have:
\begin{equation} \lb{hp}
|h(t)|^p = \s'(t) \ \ \ a.e.
\end{equation}

Suppose that $G_\f\supsetneq U.$ Then $G_\f \subseteq NS(\bbR_+,a)$ for some $a>1$, i.e. there exists $b>0$ with $\s'(t)\in\{ba^k\}^{k\in\bbZ}$ for all $t$.

Consider $\s\in G_\f\setminus U.$ \Buo we can assume that $\m(\{t:\ \s'(t) = b\})>0$. Then by \eqref{hp} and \eqref{fi} we get
\begin{equation*}
\f(b^{1/p}x) = b\f(x) \text{\ \ \ \ \ for all $x$}.
\end{equation*}
In particular, when $x=1$, we get $\f(b^{1/p}) = b.$ Therefore
\begin{equation} \lb{bp}
\f(b^{1/p}x) = \f(b^{1/p})\f(x) \text{\ \ \ \ \ for all $x$}.
\end{equation}

Hence $b^{1/p}$ belongs to the multiplier group of $\f.$
If  the multiplier group of $\f$ is trivial (i.e. $ = \{1\}$) then $b=1$ for all $\s\in G_\f\setminus U,$ which is impossible.
Thus the multiplier group of $\f$ is cyclic, generated, say, by $\bar a$.

Let $k\in\bbZ$ be such that $\m(\{t:\ \s'(t) = ba^k\})>0$. Similarly as above we get by \eqref{hp}, \eqref{fi} and \eqref{bp}
\begin{align*}
\f(b^{1/p}a^{k/p}x) &= ba^k\f(x) = a^k\f(b^{1/p}x)\\
                    &= \f(a^{k/p})\f(b^{1/p}x)
\end{align*}
for all $x\ge0$.

Thus $a^{k/p}$ belongs to the multiplier group of $\f,$ i.e. $a^{k/p} = \bar a^n$ for some $n\in\bbN.$ Hence $G_\f = NS(\bar a^p,1)$.
\end{pf*}

\begin{pf*}{Proof of Corollary~\ref{orlicz}}
Assume first that $L_{\varphi}$ and $L_{\psi}$ are isometric.

Since $0<\varphi(t)<\infty$ and $0<\psi(t)<\infty$ for all $t>0 $, $L_\varphi$ and $L_\psi$  both satisfy $(GP)$. Thus, by  Corollary~\ref{az}, ${\Iso}(L_\varphi) = \Iso(L_\psi).$ 

Further, if  ${\Iso}(L_\varphi) = \Iso(L_\psi) \neq NS(a,1)$ for any $a>1,$ then the identity map  between $L_{\varphi}$ and $L_{\psi}$ is also an isometry and thus, by Theorem~\ref{LO}, we have $(1)$.

Thus suppose that  ${\Iso}(L_\varphi) = \Iso(L_\psi) = NS(a,1)$ for some $a>1.$  Then, by Lemma~\ref{gror}, there exists $p$ ($1\le p<\infty$) so that $L_{\varphi}$, $L_{\psi}$ are isomorphic to $L_p$ and $\bar a = a^{1/p}$ is a generator of the multiplier group for $\psi$ (resp. $\f$) i.e.
$$ \psi(\bar at) = \psi(\bar a)\psi(t),    \text{\ \ \ \ \  for all $t$}$$
and if  $ \psi(ct) = \psi(c)\psi(t)$ for all $t$, then $c=\bar a^n$ for  some $n\in\bbZ.$

Now, since $L_\psi$ is isomorphic to $L_p$, we can represent $\psi$ as:
$$\psi(t) = t^p\cdot \tilde\psi(t) \text{\ \ \ \ \ for all $t$},$$
where, by \cite[Theorem~8.1]{Kr-Rut} 
\begin{equation} \lb{bound}
C_1\le\tilde\psi(t) \le C_2 \text{\ \ \ \ \ for all $t\ge t_0$}, 
\end{equation}
for some $C_1, C_2, t_0 >0$ (and $\tilde\psi(t) \not\equiv const.$ since $L_\psi$ is not isometric to $L_p$).

It is easy to see that 
$$\tilde\psi(ct) = \tilde\psi(c)\tilde\psi(t)  \iff  \psi(ct) = \psi(c)\psi(t).$$
Thus $\bar a$ is a generator of the multiplier group for $\tilde\psi.$ Further
$$ \tilde\psi(\bar a^n) = (\tilde\psi(\bar a))^n \text{\ \ \ \ \ for all $n\in\bbN.$}$$
Therefore, by \eqref{bound}, we get $\tilde\psi(\bar a) = 1.$

Now suppose that  $T:L_\f \lra L_\psi$ is a surjective isometry. Then by Theorems~\ref{main} and \ref{Z}(c) there exist $\s\in NS(\bbR_+,a)$ and a map $h$ on $[0,1]$ with $|h(t)|^p = \s'(t) \ \ a.e.$, such that
$$Tf = h\cdot f\circ \s  \text{\ \ \ \ \ for all $f\in L_\f.$}$$

Thus
for all $f\in L_\f$ we have:
\begin{align*}
\|f\|_\f &= \inf\left\{c\ : \int_0^1 \psi\left( \frac{|h(t)| \ |f(\s(t))|}{c}\right) dt \le 1 \right\} \\
                &= \inf\left\{c\ : \int_0^1 \frac{|h(t)|^p \ |f(\s(t))|^p}{c^p}\tilde\psi\left( \frac{|h(t)| \ |f(\s(t))|}{c}\right) dt \le 1 \right\} \\
\intertext{letting $s = \s(t)$, since  $|h(t)|^p = \s'(t)$, we get }
                 &= \inf\left\{c\ : \int_0^1 \frac{ |f(s)|^p}{c^p}\tilde\psi\left( \frac{|h(\s^{-1}(s))| \ |f(s)|}{c}\right) ds \le 1 \right\} \\
                 &= \inf\left\{c\ : \int_0^1 \frac{ |f(s)|^p}{c^p}\tilde\psi\left( \frac{b^{1/p} a^{k(s)/p} \ |f(s)|}{c}\right) ds \le 1 \right\} \\
\intertext{where $k(s) \in\bbZ$ is such that $\s'(\s^{-1}(s)) = b a^{k(s)}$ (since $\s\in NS(\bbR_+,a)$)}
                &= \inf\left\{c\ : \int_0^1 \frac{ |f(s)|^p}{c^p}\tilde\psi\left( \frac{b^{1/p}\ |f(s)|}{c}\right) ds \le 1 \right\} \\
\intertext{since $a^{1/p} = \bar a$ is a  multiplier  for $\tilde\psi$ and $\tilde\psi(a^{1/p}) = 1,$}
                  &= \inf\left\{c\ : \int_0^1 \psi_1\left( \frac{|f(s)|}{c}\right) ds \le 1 \right\} \\
                   &=  \|f\|_{\psi_1},
\end{align*}
where $\psi_1(x) = (1/b) \psi(b^{1/p}x),$ for all $x.$

Thus $Id: L_\f \lra L_{\psi_1}$ is an isometry and by Theorem~\ref{LO}, $\f(x) = \psi_1(x)\ \  a.e.$, i.e. $\f(x) = (1/b) \psi(b^{1/p}x).$
\end{pf*}

\section{Proof of Theorem~\ref{main}} \lb{last}
The proof of Theorem~\ref{main} follows the same essential steps as 
that of    \cite[Theorem~6.4]{KR}. 

Notice that, similarly as in \cite{KR}, it is enough to assume that spaces $X, \ Y$ satisfy property $(GP).$ Otherwise $X', Y'$ satisfy $(GP)$ and we prove that $T':Y'\lra X'$ is elementary. Then by \cite[Lemma~2.4]{KR} we get that also $T$ is elementary (see also final remarks in the proof of Theorem~6.4 in \cite[p. 322]{KR}).

We will need an analogue of    \cite[Proposition~6.3]{KR}. 
The proof is the same so we omit it. 

\begin{prop} \lb{6.3'}Let $X$ be an r.i. space on
$[0,1]$ with property $(GP)$, and such that $X\neq L_2$
(isometrically) and let $Y$ be any K\" othe function space on
$[0,1]$ for which $Y'$ is norming.
  Suppose $T:X\lra Y$ is a surjective isometry.  Then
  there exists a sequence
of Borel maps $\sigma_n:[0,1]\lra [0,1]$ and Borel functions
$a_n$ on $[0,1]$ so that $|a_n(s)|\ge |a_{n+1}(s)|$ a.e. for
every $n$ and $\sigma_m(s)\neq \sigma_n(s)$ whenever $m\neq
n$ and $s\in [0,1]$ and for which
$$Tf(s) = \sum_{n=1}^{\infty} a_n(s) f(\s_n(s)) \ \ a.e.\ \
\text{for all $f \in L_{\infty}.$}$$
and the above random measure representation (also called
abstract kernel) is unique.

Further, there exists a constant $C$, independent of the isometry
$T,$ such that
$$\int_0^1 \left( \sum_{n=1}^{\infty} |a_n(s)| \right) ds \leq C.$$
 \end{prop}

\begin{rem} 
The reader might have noticed a slight difference in the assumptions 
of Proposition~\ref{6.3'} here and Proposition~6.3 in \cite{KR}. 
In fact property $(GP)$ described here is weaker than property $(P)$ 
from \cite{KR}. 
Because of the excessive technical terminology we do not
 wish to fully 
analyze the slight differences that this produces --- we just want to 
point out that property $(P)$ was introduced to ensure the validity of 
   \cite[Lemma~5.3]{KR}, which provides the crucial step for the argument 
in \cite{KR}. 
The generalized property $(GP)$ even though weaker than property $(P)$ 
still is sufficient for \cite[Lemma~5.3]{KR} to hold (with identical proof) 
and hence there is no difference in any subsequent arguments including 
the proof of the proposition quoted above.
\end{rem} 

Next we need to show that $|a_2(s)| = 0 \ a.e.$
We follow the construction very similar to that presented in
\cite{KR}. We will recall all the necessary notation.

First, by Proposition~\ref{6.3'}, $T^{-1}$ has a random measure
representation with all the properties described above.
$$T^{-1}f(s) = \sum_{n=1}^{\infty} b_n(s) f(\r_n(s)) \ \ a.e.\ \
\text{for all $f \in L_{\infty},$}$$
and there exists $K>0,$  depending only on $Y,$ such that
$$\int_0^1 \left( \sum_{n=1}^{\infty} |b_n(s)| \right) ds \leq K.$$

Let $M_N(s)$ be the
greatest index such that $\sigma_1(s),\ldots,\sigma_M(s)$
belong to distinct dyadic intervals $D(N,k).$ Then
$M_N(s)\lra\infty$ for all $s$ and it follows easily that
given $\e>0$ we can find $M,N$ and a Borel subset $E$
of $[0,1]$ with $\lambda(E)>1-\e$ and such that
$M_N(s)\ge M$ for $s\in E,$ and
\begin{align*}
\int_{[0,1]\setminus E}\sum_{n=1}^{\infty}|a_n(t)|dt &<\e\\
\int_E\sum_{n=M+1}^{\infty}|a_n(s)| ds&<\e.
\end{align*}

Next we use the group $\cal T$ defined in \cite{KR}. We
recall the definition.

Set $P=2^N$ and let us
identify the circle group $\bbT$ with $\bbR / \bbZ
=[0,1)$ in the natural way.  For $\theta\in [0,1)^{P}$ we
define a measure preserving Borel automorphism
$\gamma=\gamma_1,\ldots,\theta_P)$ given by
$\gamma(0)=0$ and then $$
\gamma(s)=s+(\theta_k-\rho)2^{-N}$$ for $(k-1)2^{-N}<s\le
k2^{-N}$ where $\rho=1$ if $2^Ns+\theta_k>k$ and $\rho=0$
otherwise.  Thus $\gamma$ leaves each $D(N,k)$ invariant.
The set of all such $\gamma$ is a group of automorphisms
$\Gamma$ which we endow with the structure of the
topological group $\bbT^P=[0,1)^P.$ We denote Haar
measure on $\Gamma$ by $d\gamma.$ For each $k$ let
$\Gamma_k$ be the subgroup of all $\gamma(\theta)$ for which
$\theta_i=0$ when $i\neq k.$ Thus
$\Gamma=\Gamma_1\ldots\Gamma_P.$

We also let the finite permutation group $\Pi_P$ act on
$[0,1]$ by considering a permutation $\pi$ as inducing an
automorphism also denoted $\pi$ by $\pi(0)=0$ and then
$\pi(s)=\pi(k)-k+s$ for $(k-1)2^{-N}<s\le k2^{-N}.$ We again
denote normalized Haar measure on $\Pi_P$ by $d\pi.$ Finally
note that the set $\Gamma\Pi_P=\cal T$ also forms a compact
group when we endow this with the product topology and Haar
measure $d\tau=d\gamma\,d\pi$ when $\tau=\gamma\pi.$

We now wish to consider the isometries $V_{\tau}:X\lra X$ for
$\tau\in\cal T$ defined by $V_{\tau}x=x\circ\tau.$ For every
$\tau\in\cal T$ the operator $S(\tau)=TV_{\tau}T^{-1}$ is a
surjective isometry of $Y$ and so, by   \cite[Theorem~6.4]{KR},
the unique random measure representation of $S(\tau)$ consists
of exactly one term i.e. there exist Borel function $h$ on $[0,1]$
and an invertible Borel map $\s:[0,1] \lra [0,1]$ such that
$S(\tau)f = hf\circ\s$ for all $f$ in $Y.$

On the other hand the abstract kernels of $T^{-1},\; V_{\tau}$ and
$T$ `multipy' i.e., more precisely, the following version of
 \cite[Lemma~6.5]{KR} holds (with the identical proof).

\begin{lem} \lb{6.5'}For almost every $\tau\in \cal T$ we
have that  \begin{align}
\int_0^1\sum_{n=1}^{\infty}\sum_{j=1}^{\infty}|a_j(s)|
|b_n(\tau\sigma_js)|
ds &<\infty \lb{(1)}\\
\sum_{n=1}^{\infty}\sum_{j=1}^{\infty}a_j(s)b_n(\tau\sigma_js)
\delta(
\r_n\tau \sigma_js) &= \nu^{S(\tau)}_s\text{ a.e.}\lb{(2)}
\end{align}  \end{lem}

Finally, we state an equivalent of  \cite[Lemma~6.6]{KR} in the
formulation more convenient for our present purpose.

\begin{lem} \lb{6.6'} For almost every $(s,\tau) \in
E\times\cal T$ and for any two distinct pairs $(n,j),(m,i)$ where
$m,n\in\bbN$ and $i,j\le M$ if
$\r_n\tau\sigma_js=\r_m\tau\sigma_is$ then
$a_j(s)b_n(\tau\sigma_js)= 0$ or
$a_i(s)b_m(\tau\sigma_is)= 0.$
 \end{lem}


We now complete the proof that $|a_2(s)| = 0 \ a.e.$

By Lemma~\ref{6.5'} we have that for almost every $(s,\tau) \in
E\times\cal T$ the series
$$ \sum_{n=1}^{\infty}\sum_{j=1}^{\infty}a_j(s)b_n(\tau\sigma_js)
\delta(\r_n\tau \sigma_js)$$
has lenght one.
 Hence for all, except one pair $(j,n)\in \bbN\times\bbN$ we have
\begin{equation}
\label{3}
a_j(s)b_n(\tau\sigma_js)
=
-\sum_{k,l\in
I(j,n,\tau,s)}a_k(s)b_l(\tau\sigma_ks),
\end{equation}
where $I(j,n,\tau,s)=\{(k,l)\,:\,\r_n\tau\sigma_js = \r_l\tau\sigma_ks \}.$

We restrict our attention to $j = 1,2$ and $n = 1.$ Then \eqref{3}
is valid for at least one of the pairs (1,1), (2,1).
By Lemma~\ref{6.6'}  for almost every $(s,\tau) \in
E\times\cal T$ if $(k,l) \in I(j,1,\tau,s)$ then $k> M.$
Hence, if \eqref{3} holds for (j,1) we have
\begin{equation*}
\begin{split}
\int_{E\times\cal T} |a_j(s)b_1(\tau\sigma_js)| ds d\tau &=
\int_{E\times\cal T} |\sum_{k,l \in I(j,n,\tau,s)}
a_k(s)b_l(\tau\sigma_ks)| ds d\tau  \\
&\le \int_{E\times\cal T} \sum_{k>M}\sum_{l=1}^{\infty}
|a_k(s)||b_l(\tau\sigma_ks)| ds d\tau \\
&=  \int_E \left(\sum_{k>M}|a_k(s)| \int_{\cal T}\sum_{l=1}^{\infty}
|b_l(\tau\sigma_ks)| d\tau\right)  ds
\end{split}
\end{equation*}
But for almost every $s$  
${\disp  \int_{\cal T} |b_l(\tau\sigma_ks)| d\tau = \int_0^1
|b_l(t)| dt.}$ Thus we obtain:
\begin{equation*}
\int_{E\times\cal T} |a_j(s)b_1(\tau\sigma_js)| ds d\tau \le
\int_E  \sum_{k>M}|a_k(s)| ds \int_0^1\sum_{l=1}^{\infty}
|b_l(t)| dt   \le \e K.
\end{equation*}

On the other hand:
\begin{equation*}
\begin{split}
\int_{E\times\cal T} |a_j(s)&b_1(\tau\sigma_js)| ds d\tau
= \int_E   |a_j(s)| \int_{\cal T} |b_1(\tau\sigma_js)| d\tau ds \\
&= \int_E   |a_j(s)| ds \int_0^1 |b_1(t)| dt  = \be \int_E   |a_j(s)|
ds,
\end{split}
\end{equation*}
where ${\disp \be = \int_0^1 |b_1(t)| dt }$ ($0<\be \le K$).

Therefore
$$ \int_E   |a_j(s)| ds \le \frac{\e K}{\be}$$
and
$$ \int_0^1   |a_j(s)| ds = \int_E   |a_j(s)| ds +
\int_{[0,1]\setminus E}  |a_j(s)| ds \le \frac{\e K}{\be} + \e.$$

Since $\e >0$ was chosen arbitrarily we see that \eqref{3}
implies that ${\disp \int_0^1   |a_j(s)| ds = 0.}$ Therefore
\eqref{3} cannot hold for (1,1) and, thus, it holds for (2,1)
and $|a_2(s)| = 0 \ \ a.e. \qed $

  \bibliographystyle{plain}
  \bibliography{tref}

  \end{document}